\documentclass{amsart}
 
\RequirePackage{amsmath}
\RequirePackage{amscd}
\RequirePackage{amssymb}
\RequirePackage{xy}
 
\xyoption{all}
\CompileMatrices
 
\theoremstyle{plain}
\newtheorem{theorem}{Theorem}[section]
\newtheorem{proposition}[theorem]{Proposition}
\newtheorem{lemma}[theorem]{Lemma}
\newtheorem{corollary}[theorem]{Corollary}

\newtheorem*{mainthm}{Theorem}
\newtheorem*{question}{Question}
 
\theoremstyle{remark}

\numberwithin{equation}{section}

\newcommand{\Q}{{\mathbf Q}}
\newcommand{\Qp}{\Q_{p}}
\newcommand{\Z}{{\mathbf Z}}
\newcommand{\Zp}{\Z_{p}}
\newcommand{\Zpp}{\Z_{(p)}}

\renewcommand{\a}{{\mathfrak a}}
\newcommand{\g}{{\mathfrak g}}

\newcommand{\p}{{\mathfrak p}}
\newcommand{\pt}{\tilde{\mathfrak p}}

\renewcommand{\O}{{\mathcal O}}
\newcommand{\Ot}{\tilde{\O}}

\newcommand{\Fa}{F_{\alpha}}
\newcommand{\Fbar}{\bar{F}}

\newcommand{\Poincare}{Poincar\'e}
\newcommand{\ts}{\textstyle}

\newcommand{\Ga}{G_{\alpha}}
\newcommand{\Gpi}{G_{p^{\infty}}}

\newcommand{\ga}{\g_{\alpha}}
\newcommand{\gpi}{\g_{p^{\infty}}}
\newcommand{\gpn}{\g_{p^{n}}}
\newcommand{\gas}{\ga(\Gamma)}
\newcommand{\gpis}{\gpi(\Gamma)}
\newcommand{\gpns}{\gpn(\Gamma)}

\newcommand{\la}{\lambda_{\alpha}}
\newcommand{\lpi}{\lambda_{p^{\infty}}}
\newcommand{\lpn}{\lambda_{p^{n}}}
\newcommand{\las}{\la^{\Gamma}}
\newcommand{\lax}{\la^{\O \cdot x}}
\newcommand{\lpis}{\lpi^{\Gamma}}
\newcommand{\lpns}{\lpn^{\Gamma}}

\newcommand{\soa}{\ts {\frac{1}{\alpha}\Gamma}}
\newcommand{\xoa}{{\ts \frac{x}{\alpha}}}

\newcommand{\tors}{\mathrm{tors}}

\newcommand{\cd}{\!\cdot\!}
\newcommand{\inj}{\hookrightarrow}
\newcommand{\invlim}{\varprojlim}
\newcommand{\op}[2]{\underset{#1}{\overset{#2}{\oplus}}}
\newcommand{\ST}{\SelectTips{cm}{}}
\newcommand{\surj}{\twoheadrightarrow}
\renewcommand{\th}{\text{th}}
\newcommand{\too}{\longrightarrow}

\DeclareMathOperator{\coker}{coker}
\DeclareMathOperator{\End}{End}

\DeclareMathOperator{\Frob}{Frob}
\DeclareMathOperator{\Gal}{Gal}

\DeclareMathOperator{\Hom}{Hom}
\DeclareMathOperator{\im}{im}
\DeclareMathOperator{\ord}{ord}
\DeclareMathOperator{\rank}{rank}
\DeclareMathOperator{\red}{red}
\DeclareMathOperator{\res}{res}
\DeclareMathOperator{\Tor}{Tor}

\newcommand{\Ef}{\End_{F}\!}
\newcommand{\rv}{\red_{v}\!}
\newcommand{\rw}{\red_{w}\!}

\title[Kummer theory and reductions of Mordell-Weil groups]
{Kummer theory of abelian varieties and reductions of Mordell-Weil groups}
\author{Tom Weston}
\email{weston@math.berkeley.edu}

\begin{document}

\maketitle

Let $A$ be an abelian variety over a number field $F$.
We write 
$\rv : A(F) \to A(k_{v})$ for the reduction map
at a place $v$ of $F$ with residue field $k_{v}$.
W.\ Gajda has posed the following question.

\begin{question}
Let $\Sigma$ be a subgroup of $A(F)$.
Suppose that $x$ is a point of
$A(F)$ such that $\rv x$ lies in $\rv \Sigma$
for almost all places $v$ of $F$.  Does it then follow that $x$ lies in
$\Sigma$?
\end{question}

In this paper we use methods of Kummer theory to provide the following
partial answer to this question.

\begin{mainthm}
Let $A$ be an abelian variety over a number field $F$ and assume
that $\Ef A$ is commutative.  Let $\Sigma$ be a subgroup of $A(F)$ and
suppose that $x \in A(F)$ is such that
$\rv x \in \rv\Sigma$
for almost all places $v$ of $F$.  Then $x \in \Sigma + A(F)_{\tors}$.
\end{mainthm}

It does not appear that the torsion ambiguity can be eliminated with our
present approach, and it is not clear to the author
how to modify the arguments for the non-commutative case.
We note that our theorem applies in particular to products of non-isogenous
elliptic curves.

Gajda's question has its origins in the support problem of P.\ Erd\H{o}s:
if $x$ and $y$ are positive integers such that for any $n \geq 1$
the set of primes dividing $x^{n}-1$ is the same as the set of primes dividing
$y^{n}-1$, then must $x$ equal $y$? 
Corrales-Rodrig\'a\~nez and Schoof
gave an affirmative answer to this question in \cite{CRS} and also answered
the corresponding question for elliptic curves;
this was generalized by Banaszak, Gajda and Kraso\'n in \cite{BGK}
to certain abelian varieties with complex or real
multiplication and $\Ef A$ a commutative maximal order.
In this context the support problem takes the following form.

\begin{question}
Let $x,y \in A(F)$ be non-torsion points.
Suppose that the order of $\rv x$ divides the order of
$\rv y$ for almost all places $v$ of $F$.
Does it follow that $x$ and $y$ satisfy an $\Ef A$-linear relation
in $A(F)$?
\end{question}

Taking $\Sigma = \Ef A \cd y$,
the support problem implies a weak form of our main theorem in the
case that $\Sigma$ is a cyclic $\Ef A$-module.
The more precise question of Gajda we consider is
one possible modification of the support problem
for abelian varieties to a non-cyclic setting.
The approach we use here is quite different from that of
\cite{CRS} and \cite{BGK},
relying more on the study of the Mordell-Weil group of $A$ as
a module for $\Ef A$ and less on Galois cohomology.

We give now an overview of our argument in the simplest case.
Assume that $A$ is simple, that $\O := \Ef A$ is integrally closed (so
that it is a Dedekind domain), and that $A(F)$ is a free $\O$-module.
With $\Sigma \subseteq A(F)$ and $x \in A(F)$ as in the theorem, it suffices
to prove that $x \in \Sigma \otimes \Zpp$ for every prime $p$ (with $\Zpp$
the localization of $\Z$ away from $p$).
Fix, then, a prime $p$ and suppose that $x \notin \Sigma \otimes \Zpp$.
The first step, which is purely algebraic, is to show that under this
assumption one can choose an $\O$-basis $y_{1},\ldots,y_{r}$ of $A(F)$
such that $\psi_{1}(x) \notin \psi_{1}(\Sigma) + p^{a}\O$ for some $a > 0$;
here $\psi_{1} : A(F) \to \O$ is the projection onto the $y_{1}$-coordinate.

The next step is to choose an appropriate place $v$ of $F$.
We work instead over the extensions $F(A[p^{n}])$ of $F$.
Using Kummer theory and the Cebatorev density theorem, we show that
there is a $b > 0$ such that for any sufficiently large
$n$ there is a place $w$ of $F(A[p^{n}])$ with
$\rw  y_{2},\ldots,\rw y_{r} \in p^{n}A(k_{w})$, while
$\rw  y_{1} \notin \p_{i}^{b}A(k_{w})$ for any $i$; here
$p\O = \p_{1}^{e_{1}}\cdots\p_{g}^{e_{g}}$ is the ideal factorization of
$p$ in $\O$.

Fix $n \geq a+b$ and choose such a place $w$.
By hypothesis we have $\rw  x = \rw y$ for some $y \in \Sigma$.
Expanding in terms of our chosen basis of $A(F)$, the choice of $w$
implies that
$$\bigl(\psi_{1}(x)-\psi_{1}(y)\bigr)\rw y_{1} \in p^{n}A(k_{w}).$$
On the other hand, using the properties of $\psi_{1}$ and of $w$, one
can show directly that
$$(\psi_{1}(x)-\psi_{1}(y))\rw y_{1} \notin p^{a+b}A(k_{w}).$$
As $n \geq a+b$, we have a contradiction,
so that we must have had
$x \in \Sigma \otimes \Zpp$.  This completes our sketch of the argument
in this case.

We now review the contents of this paper in more detail.
We begin in Section~\ref{s:k1} with a review of Kummer theory and in
Section~\ref{s:k2} we
adapt the methods of Bashmakov--Ribet as in \cite{Ribet} to prove that
the cokernel of the $p$-adic Kummer map is bounded.  In Section~\ref{s:k3}
we discuss the relation between Kummer theory and reduction maps.

In the sketch above we assumed that $\O$ was an integrally closed domain
and that $A(F)$ was free over $\O$.  The algebra required to eliminate
these assumptions is developed in Section~\ref{s:a}.
These results are combined with Kummer theory to produce places $w$
as above in Section~\ref{s:a3}, and the proof of our main theorem is
given in Section~\ref{s:a4}.

The author wishes to thank Ken Ribet for suggesting this problem and Mark
Dickinson for helpful conversations.  The author would also like to 
acknowledge
the \texttt{gp-pari} computing package, as it was in the course of some 
seemingly unrelated calculations 
that some of the main ideas of this paper became apparent.

\section{Kummer theory} \label{s:k}

\subsection{Review of Kummer theory} \label{s:k1}

Let $A$ be an abelian variety over a number field $F$;
set $\O = \Ef A$.  For $\alpha \in \O$ we set
$\Fa = F(A[\alpha])$ and $\Ga = \Gal(\Fa/F)$.
The {\it Kummer map} 
$$\kappa_{\alpha} : A(F)/\alpha \to \Hom_{\Ga}\bigl(\Gal(\Fbar/\Fa),
A[\alpha]\bigr)$$
is defined as the composition
$$A(F)/\alpha \inj H^{1}\bigl(F,A[\alpha]\bigr) \overset{\res}{\too}
H^{1}\bigl(\Fa,A[\alpha]\bigr)^{\Ga}$$
with the first map a coboundary map for the $\Gal(\Fbar/F)$-cohomology of
the Kummer sequence
$$0 \to A[\alpha] \to A(\Fbar) \overset{\alpha}{\too} A(\Fbar)
\to 0$$
and the second map restriction to $\Fa$.
(Concretely, for $x \in A(F)$, $\kappa_{\alpha}(x)$ is the
homomorphism sending $\gamma \in \Gal(\Fbar/\Fa)$ to $\gamma(\xoa)
- \xoa \in A[\alpha]$ where $\xoa$ is some fixed $\alpha^{\th}$-root of
$x$ in $A(\Fbar)$.)

If $\Gamma$ is an $\O$-submodule of $A(F)$ and $\alpha \in \O$, we
write $\Fa(\soa)$ for the extension of $\Fa$
generated by all $\alpha^{\th}$-roots of elements of $\Gamma$;
alternately, $\Fa(\soa)$ is the fixed field of
the intersection of the kernels of the homomorphisms $\kappa_{\alpha}(\Gamma)$.
The Galois group $\gas := \Gal(\Fa(\soa)/\Fa)$ is an $\O[\Ga]$-module and
$\kappa_{\alpha}$ restricts to an $\O$-linear map
$\Gamma \to \Hom_{\Ga}\bigl(\gas,A[\alpha]\bigr)$.
We write the $\O[\Ga]$-dual of this map as
$$\las : \gas \inj \Hom_{\O}\bigl(\Gamma,A[\alpha]\bigr).$$

\subsection{$p$-adic Kummer theory} \label{s:k2}

Fix a rational prime $p$; set
$\O_{p} = \O \otimes \Zp$ and $K_{p} = \O \otimes \Qp$.
The Tate module $T_{p}A := \invlim A[p^{n}]$ 
(resp.\ Tate space $V_{p}A := T_{p}A \otimes_{\Zp} \Qp$) is
naturally an $\O_{p}[\Gpi]$-module 
(resp.\ $K_{p}[\Gpi]$-module) where $\Gpi = \Gal
\bigl(F(A[p^{\infty}])/F\bigr)$.
It follows from \cite[Section 19, Corollary 2]{Mumford} that there is
a decomposition
\begin{equation} \label{eq:csa}
K_{p} = {\ts \prod} M_{n_{i}}K_{i}
\end{equation}
where $M_{n_{i}}K_{i}$ is the central simple algebra of 
$n_{i} \times n_{i}$-matrices over the division ring $K_{i}$.
Corresponding to (\ref{eq:csa}) is a decomposition
$V_{p}A = \oplus V_{i}A^{n_{i}}$
of $V_{p}A$ into $K_{i}[\Gpi]$-modules.  By 
\cite[Theorem 4]{Faltings}, we have 
\begin{equation} \label{eq:end}
\End_{\Qp[\Gpi]}V_{i}A = K_{i}
\end{equation}
for each
$i$; in particular, each $V_{i}A$ is an irreducible $K_{i}[\Gpi]$-module.
We record a second immediate consequence of (\ref{eq:end}) in the next lemma.

\begin{lemma} \label{l:fin3}
Let $\Gamma$ be an $\O$-module.  Then the evaluation map
$$\Gamma \otimes_{\O} 
K_{i} \to \Hom_{K_{i}[\Gpi]}\bigl(\Hom_{\O}(\Gamma,V_{i}A),
V_{i}A\bigr)$$
is an isomorphism.
\end{lemma}

Fix an $\O$-submodule $\Gamma$ of $A(F)$.
The inverse limit $\gpis$ of the $\gpns$ is naturally an $\O_{p}[\Gpi]$-module
endowed with an injection
$$\lpis : \gpis \inj \Hom_{\O}(\Gamma,T_{p}A).$$
More generally, since $\O_{p}/p^{n} \cong \O/p^{n}$ for all $n$, for any
$\O$-module $\Gamma \subseteq A(F) \otimes \Zp$ we can still
define $\gpns$ and $\lpns$ for $n \leq \infty$.
In any case, there is a $K_{p}[\Gpi]$-module decomposition
\begin{equation} \label{eq:gi}
\gpis \otimes_{\Zp} \Qp = \oplus \g_{i}(\Gamma)^{n_{i}}
\end{equation}
(with $n_{i}$ as in (\ref{eq:csa}))
into $K_{i}[\Gpi]$-modules, and there are natural injections
$$\lambda_{i}^{\Gamma} : \g_{i}(\Gamma) \inj \Hom_{\O}(\Gamma,V_{i}A).$$
The decomposition (\ref{eq:gi}) is
functorial in the sense that there is a natural surjection
$\g_{i}(\Gamma) \surj \g_{i}(\Gamma')$ for any $\O$-submodule $\Gamma'$ of
$\Gamma$.

The main result of Kummer theory we need is the following.
The proof is a straightforward adaptation of
the methods of Bashmakov and Ribet.

\begin{proposition} \label{p:padic}
Fix a rational prime $p$ and
let $\Gamma$ be an $\O$-submodule of $A(F)$.
Then the cokernel of $\lpns$ is bounded independent of $n$.
\end{proposition}
\begin{proof}
First consider the cyclic case $\Gamma = \O \cd x$ for $x \in A(F)$.
If $\Gamma \cong \O$, then $\Z \cd x$ 
is Zariski dense in $A$; the proposition
thus follows from
\cite[Theorem 2]{Bertrand} in this case.  More generally, let $A'$ denote the
largest abelian subvariety of $A$, defined over $F$, in which $\Z \cd x$ is
Zariski
dense; set $\O' = \Ef A'$.  Using the \Poincare~ reducibility theorem
(see \cite[Section 19, Theorem 1]{Mumford}), one checks easily that
$$\Hom_{\O}(\Gamma,V_{p}A) \cong \Hom_{\O'}(\Gamma,V_{p}A'),$$
so that the general cyclic  case follows from \cite[Theorem 2]{Bertrand}
applied to $A'$.
In fact, one has $\coker \lpi^{\O \cd x} = \coker \lpi^{\O \cd x'}$
whenever $x,x' \in A(F)$ are sufficiently $p$-adically congruent, so that
the same arguments apply for arbitrary $x \in A(F) \otimes \Zp$.

For general $\Gamma$
it suffices to show that each of the injections $\lambda_{i}^{\Gamma}$
is an isomorphism.
Suppose, then, that some $\lambda_{i}^{\Gamma}$ is not surjective.
Since $\Hom_{\O}(\Gamma,V_{i}A)$ is a direct sum of copies of the
irreducible $K_{i}[\Gpi]$-module $V_{i}A$
(and thus in particular is a semisimple $K_{i}[\Gpi]$-module), it follows
that there exists a $K_{i}[\Gpi]$-surjection
$$\varphi : \Hom_{\O}(\Gamma,V_{i}A) \surj V_{i}A$$
annihilating $\g_{i}(\Gamma)$.
By Lemma~\ref{l:fin3} the map $\varphi$ is given by evaluation
at some $x \in \Gamma \otimes_{\O} K_{i}$; using the injection $K_{i} \inj
K_{p}$ and scaling $\varphi$ if necessary,
we may in fact assume that $x \in \Gamma \otimes \Zp$.
There is then a commutative diagram
$$\ST\xymatrix{
{\g_{i}(\Gamma)} 
\ar@{^{(}->}[r]^-{\lambda_{i}^{\Gamma}} 
\ar@{->>}[d] & {\Hom_{\O}(\Gamma,V_{i}A)}
\ar@{->>}[d]^{\varphi} \\
{\g_{i}(\O \cd x)} 
\ar@{^{(}->}[r]^-{\lambda_{i}^{\O \cd x}} & {V_{i}A}}$$
The clockwise composition is zero by construction, so
that we must have $\lambda_{i}^{\O \cd x} = 0$ as well.
By the cyclic case considered above this implies that $x$ maps to zero
in $\Gamma \otimes_{\O} K_{i}$.
But then $\varphi$, which is evaluation at $x$, is also zero.  This contradicts
the surjectivity of $\varphi$ and thus proves the proposition.
\end{proof}

\subsection{Reductions and Frobenius elements} \label{s:k3}

We write $k_{w}$ for the residue field of a finite extension $F'$ of $F$
at a place $w$ and
$\rw  : A(F') \to A(k_{w})$ for the reduction map.

\begin{lemma} \label{l:frob}
Fix $\alpha \in \O$ and $x \in A(F)$.
Let $w$ be a finite place of
$\Fa$, relatively prime to $\alpha$, at which $A$ has good reduction.
Then $\rw x$ lies in $\alpha A(k_{w})$ if and only if
$\lax(\Frob_{w}) = 0$, where $\Frob_{w} \in \Gal(\Fa(\xoa)/\Fa)$
is the Frobenius element at $w$.
\end{lemma}
\begin{proof}
Fix an $\alpha^{\th}$-root $\xoa$ of $x$ in $A(\Fbar)$ and a place $w'$ of
$\Fa(\xoa)$ over $w$.  If $\lax(\Frob_{w})=0$, then
$w'$ is completely split over $w$ so that
$k_{w'} = k_{w}$.  In particular, $\red_{w'}\!\xoa \in A(k_{w'})$ lies in
$A(k_{w})$; thus $\rw x \in \alpha A(k_{w})$ as claimed.

Conversely, if there is $y \in A(k_{w})$ with $\alpha y = \rw x$,
then $y - \red_{w'}\! \xoa$ lies in $A[\alpha]$.  Since $y$ and 
$A[\alpha]$ are both in $A(k_{w})$ we conclude that $\red_{w'}\!\xoa$ is
in $A(k_{w})$ as well.  In particular, we have
\begin{equation} \label{eq:frob}
\Frob_{w} \bigl(\red_{w'}\!\xoa\bigr) - \red_{w'}\!\xoa = 0.
\end{equation}
On the other hand,
$\Frob_{w}(\xoa) - \xoa$ already lies in $A[\alpha]$, which
injects into $A(k_{w'})$; (\ref{eq:frob}) thus forces
$$\Frob_{w} \bigl(\xoa\bigr) - \xoa = 0 \text{~in~} A(\Fbar).$$ 
This says exactly that $\lax(\Frob_{w})=0$, as claimed.
\end{proof}

We assume now that $\O$ is commutative.
Suppose that $\a$ is an ideal of $\O$ such that
$\beta \a \subseteq \alpha \O$ for some $\alpha,\beta \in \O$.
Multiplication by $\beta$ then yields a map
$A[\alpha] \to A[\a]$.

\begin{lemma} \label{l:crt2}
Let $\alpha,\beta,\a$ be as above and fix $x \in A(F)$.
Let $w$ be a finite place of 
$\Fa$, relatively prime to $\alpha$, at which $A$ has good reduction.
If $\beta\cdot\lax(\Frob_{w}) \neq 0$, then
$\rw  x \notin \a A(k_{w})$.
\end{lemma}
\begin{proof}
We prove the contrapositive.
Suppose that $\rw x \in \a A(k_{w})$.  Then
$$\beta \rw x \in \beta\a A(k_{w}) \subseteq \alpha A(k_{w}),$$
so that there is $y \in A(k_{w})$ with $\beta\rw x = \alpha y$.
On the other hand, fixing an $\alpha^{\th}$-root $\xoa$ of $x$ in
$A(\Fbar)$ and a place $w'$ of $\Fa(\xoa)$ lying above $w$,
we also have $\beta\rw x = \alpha \beta \red_{w'}\!\xoa$.
Therefore
$$y - \beta\red_{w'}\!\xoa \in A[\alpha].$$
From here the argument proceeds as in the second half of the proof of
Lemma~\ref{l:frob} above to show that $\beta\cdot\lax(\Frob_{w})=0$.
\end{proof}

We remark that the converse of Lemma~\ref{l:crt2} holds in the case that
$\alpha \O = \a \a'$ with $\a,\a'$ relatively prime and 
$\beta \in \a' \cap (1 - \a)$.

\section{Modules over commutative, reduced, finite, flat $\Z$-algebras}
\label{s:a}

\subsection{Projections} \label{s:a1}

Let $\O$ be a commutative, reduced, finite, flat $\Z$-algebra.
The normalization $\Ot$ of $\O$ decomposes as a product
$\prod_{j=1}^{h} \Ot_{j}$ of Dedekind domains.  We say that a $\Z$-linear 
map $t : \O \to \Z$ is {\it full}
if it is non-trivial on $\O \cap \Ot_{j}$ for each $j$.

\begin{lemma} \label{l:gor}
Fix a full map $t : \O \to \Z$.  Then the map
\begin{align}
\Hom_{\O}(N,\O) &\to \Hom_{\Z}(N,\Z) \label{eq:gor} \\
f &\mapsto t \circ f \notag
\end{align}
has finite cokernel for any finitely generated $\O$-module $N$.
\end{lemma}
\begin{proof}
Since $\O$ has finite index in $\Ot$, it suffices to prove the result after
replacing $\O$ by $\Ot$ and $N$ by $N \otimes_{\O} \Ot$.
We may therefore assume that $\O$
decomposes as a product $\prod \O_{i}$ of Dedekind domains.  There is then
a corresponding decomposition $N = \oplus N_{i}$, and by the definition
of a full map it suffices to prove the lemma for each
factor $N_{i}$; that is, we may assume that $\O$ is a Dedekind domain.

In this case
every finitely generated $\O$-module has a free submodule of finite index;
this allows one to reduce to the case that $N$ is free, and then to the
case that $N$ is free of rank one.  (\ref{eq:gor}) is then a map
\begin{equation} \label{eq:gor2}
\O = \Hom_{\O}(\O,\O) \to \Hom_{\Z}(\O,\Z)
\end{equation}
between two free $\Z$-modules of the same rank, so that it suffices to
prove that it is injective.  For this, note that (\ref{eq:gor2}) is
$\O$-linear; thus its kernel is an ideal of $\O$.  However, every
non-zero ideal of $\O$ has finite index and $\Hom_{\Z}(\O,\Z)$
is torsion-free; therefore (\ref{eq:gor2}) must be either zero
or injective.  As $t$ itself lies in the image, it is obviously non-zero.
\end{proof}

We now fix a finitely generated $\O$-module $N$ and a $\Z$-submodule
$M$ of $N$ containing the $\Z$-torsion submodule $N_{\tors}$ of $N$.

\begin{lemma} \label{l:notin}
Fix $x \in N$ and suppose that $p$ is a rational prime such that
$x \notin M \otimes \Zpp$.
Then there is an $\O$-linear map $\psi : N \to \O$ such
that $\psi(x) \notin \psi(M) + p^{n}\O$ for sufficiently large $n$.
\end{lemma}
\begin{proof}
Choose a $\Z$-basis $y_{1},\ldots,y_{r} \in N$ of $N/N_{\tors}$ such that
there are integers $d_{1},\ldots,d_{r}$ with
$$M = \left< d_{1}y_{1},\ldots,d_{r}y_{r} \right> \oplus N_{\tors}.$$
(Of course, some of the $d_{i}$ may be zero.)
Writing $x = a_{1}y_{1} + \cdots + a_{r}y_{r} + t$ with $a_{i} \in \Z$
and $t \in N_{\tors}$,
the fact that $x \notin M \otimes \Zpp$ 
implies that there is some index $i$ such that
\begin{equation} \label{eq:ord}
\ord_{p}a_{i} < \ord_{p}d_{i}.
\end{equation}
Let $\psi_{0} : N \to \Z$ be $\# N_{\tors}$ times projection
onto $y_{i}$; this is a well-defined map, and it follows from
(\ref{eq:ord}) that
$\psi_{0}(x) \notin \psi_{0}(M) + p^{n}\Z$ for sufficiently large $n$.
(In fact, $n > \ord_{p}  (a_{i} \cdot \# N_{\tors})$ suffices.)

Fix a full map $t : \O \to \Z$.
By Lemma~\ref{l:gor}, we can find a non-zero integer $b$
such that $b\psi_{0}$ is in the image of (\ref{eq:gor}).
Thus there is an $\O$-linear map $\psi : N \to \O$ with
$b\psi_{0} = t \circ \psi$.
Since $t(p^{n}\O) \subseteq p^{n}\Z$, we conclude that
$\psi(x) \notin \psi(M) + p^{n}\O$ for sufficiently large $n$, as desired.
\end{proof}

\subsection{Pre-bases} \label{s:a1.5}

We continue with $M \subseteq N$ as before.  Fix $y \in N$ not in
$N_{\tors}$ and
let $\varphi : \O \surj \O \cd y$ be the $\O$-linear surjection sending $1$
to $y$.
We define $\eta_{0}(y)$ to be the least positive integer $m$ such that there
exists an $\O$-linear map $\psi : \O \cd y \to \O$ with
the composition
$$\O \cd y \overset{\psi}{\too} \O \overset{\varphi}{\too} \O \cd y$$
multiplication by $m$.
(Let $K_{j}$ denote the fraction field of $\Ot_{j}$; since $\O \otimes \Q =
\prod K_{j}$, to see that
any maps $\psi$ as above exist it suffices to prove the
corresponding fact after replacing $\O$ by
$\prod K_{j}$.
In this context the map $\varphi$ idenitifes with the quotient map
$${\ts \prod} K_{j} \to \underset{j \in J}{\ts \prod} K_{j}$$
for some non-empty subset $J$ of $\{1,\ldots,h\}$,
so that the existence of $\psi$ is obvious.)

We say that $y_{1},\ldots,y_{r} \in N$
are an {\it $\O$-pre-basis} of $N$ if:
\begin{itemize}
\item $y_{i} \notin N_{\tors}$ for all $i$;
\item $(\O \cd y_{1}) \oplus \cdots \oplus (\O \cd y_{r})$
injects into $N$ with finite cokernel.
\end{itemize}
(Note that we do not require that the corresponding map
$\O^{r} \to N$ is injective.)
Let $\eta'(y_{1},\ldots,y_{r})$ be the order of this cokernel
and define
$$\eta(y_{1},\ldots,y_{r}) =
\eta'(y_{1},\ldots,y_{r})
\cdot \eta_{0}(y_{1}) \cdots \eta_{0}(y_{r}).$$
It then follows from the definition of $\eta_{0}(y_{i})$ that
there are $\O$-linear maps 
$$\psi_{i}^{y_{1},\ldots,y_{r}}: N \to \O$$
for $i=1,\ldots,r$ such that
\begin{equation} \label{eq:eta}
\eta(y_{1},\ldots,y_{r}) y = \psi_{1}^{y_{1},\ldots,y_{r}}(y)y_{1} 
+ \cdots + \psi_{r}^{y_{1},\ldots,y_{r}}(y)y_{r}
\end{equation}
for all $y \in N$.  We usually just write $\eta$ and $\psi_{i}$ if the
pre-basis $y_{1},\ldots,y_{r}$ is clear from context.
A standard inductive procedure shows that pre-bases
always exist.

\begin{proposition} \label{p:pb}
Fix $x \in N$ and suppose that
$p$ is a rational prime such that $x \notin M \otimes \Zpp$.
Then there is an $\O$-pre-basis $y_{1},\ldots,y_{r}$ of $N$ such that
$\psi_{1}(x) \notin \psi_{1}(M) + p^{n}\O$
for sufficiently large $n$.
\end{proposition}
\begin{proof}
By Lemma~\ref{l:notin}, we may choose an $\O$-linear map $\psi : N \to \O$
such that $\psi(x) \notin \psi(M) + p^{n}\O$ for sufficiently large $n$.  
Let $K'$ denote the image of $\psi \otimes \Q$; we have
$K'=\prod_{j \in J} K_{j}$ for some non-empty subset $J$ of
$\{1,\ldots,h\}$.  In particular, $K'$ is a projective $\prod K_{j}$-module,
so that there exists
a map $\varphi_{0} : K' \to N \otimes \Q$ such that
$\psi \circ \varphi_{0}$
is the identity on $K'$.  Scaling $\varphi_{0}$
by an integer we obtain an $\O$-linear map
$\varphi : \Ot' \to N$ such that $\psi \circ \varphi$ is mulitplication
by some non-zero integer; here $\Ot' = \prod_{j \in J} \Ot_{j}$.

Set $y_{1} = \varphi(1)$ and choose
an $\O$-pre-basis 
$y_{2},\ldots,y_{r}$ for $\ker \psi$.  Then $y_{1},\ldots,y_{r}$
is an $\O$-pre-basis of $N$ and 
$\psi_{1} = m\psi$ for some non-zero integer $m$.
It thus follows from the definition of $\psi$
that $\psi_{1}(x) \notin \psi_{1}(M) + p^{n}\O$
for sufficiently large $n$, as desired.
\end{proof}

\subsection{Ideals} \label{s:a2}

We continue with $\O$ as above.
Fix a rational prime $p$ and 
write the $\Z$-exponent of $\Ot/\O$ as $cp^{d}$ with $d \geq 0$ and
$c$ relatively prime to $p$.  Let
$$p\Ot = \pt_{1}^{e_{1}} \cdots \pt_{g}^{e_{g}}$$
be the factorization of $p\Ot$ into prime ideals of $\Ot$;
for each $i \in \{1,\ldots,g\}$ we let $\mu_{p}(i)$ denote the unique
$j \in \{1,\ldots,h\}$ such that 
$\pt_{i}$ is the pullback of a prime ideal on $\Ot_{j}$.
For $y \in N$ we define
$I_{p}(y) \subseteq \{1,\ldots,g\}$
to be the set of indices $i$ such that the image of $y$ in
$N \otimes_{\O} \Ot_{\pt_{i}}$ is non-torsion.
In fact, since every proper ideal of each $\Ot_{j}$ has finite index,
we have
\begin{equation} \label{eq:pbind}
I_{p}(y) = \{ i ; \rank_{\Z}\bigl((\O \cap \Ot_{\mu_{p}(i)}) \cd y\bigr)
> 0 \}.
\end{equation}

For $i=1,\ldots,g$ and any $n$, we define ideals of $\O$ by
$$\p_{i,n} = \pt_{i}^{e_{i}n} \cap \O.$$
The reader is invited to focus on the case $d = 0$, when
$\p_{i,n} = \p_{i,1}^{n}$ and the analysis below is quite a bit simpler.
In the general case, we have $cp^{d}\pt^{n}_{i} \subseteq \p_{i,n}$; since
the $\pt_{i}$ are relatively prime, it follows that
\begin{equation} \label{eq:l1}
c^{g-1}p^{d(g-1)}\O \subseteq \p_{i,n} + {\ts \underset{j \neq i}{\prod}}
\p_{j,n}
\end{equation}
for all $n$.
Furthermore, $p^{n}\Ot \cap \O \subseteq p^{n-d}\O$ for $n \geq d$, so that
\begin{equation} \label{eq:l3}
p^{n}\O \subseteq \p_{1,n} \cap \cdots \cap \p_{g,n} \subseteq p^{n-d}\O;
\end{equation}
\begin{equation} \label{eq:l2}
c^{dg}p^{n+dg}\O \subseteq \p_{1,n}\cdots\p_{g,n} \subseteq p^{n-d}\O;
\end{equation}
for any $n \geq d$.

\begin{lemma} \label{l:evil}
Let $N$ be a finitely generated $\O$-module.  Fix $\alpha \in \O$ and
$x \in N$.  Suppose that there is an index $i$ and non-negative integers
$a,b$ such that:
\begin{enumerate}
\item $\alpha \notin \p_{i,a}$;
\item $x \notin \p_{i,b}N$;
\item $N[p^{a+d}] \subseteq p^{b}N$.
\end{enumerate}
Then $\alpha x \notin p^{a+b+d}N$.
\end{lemma}
\begin{proof}
We first replace $\O$ by
$\invlim \O/\p_{i,n}$, $N$ by $\invlim N/\p_{i,n}$, and
$\Ot$ by $\invlim \O/\pt_{i}^{n}$.
Let $\pt$ denote the maximal ideal of
$\Ot$, so that $\pt^{e_{i}} = p\Ot$; set $\p_{n} = \pt^{e_{i}n} \cap \O$.  
With this notation we have
$\alpha \notin \p_{a}$ and $x \notin \p_{b}N$, and it suffices to prove
that $\alpha x \notin p^{a+b+d}N$.
Note that $\alpha \notin \pt^{e_{i}a}$, so that
there is some $\beta \in \Ot$ with $\alpha\beta = p^{a}$.

Set $C = \Ot/\O$ and $\tilde{N} = N \otimes_{\O} \Ot$; $C$ is killed
by $p^{d}$ and there is an exact sequence
\begin{equation} \label{eq:keyes}
\Tor_{1}^{\O}(N,C) \to N \overset{\iota}{\too} \tilde{N} \to 
N \otimes_{\O} C \to 0.
\end{equation}
Suppose now that $\alpha x \in p^{a+b+d}N$.  
Applying $\iota$ and multiplying by $\beta$, we find that
$p^{a}\iota(x) \in p^{a+b+d}\tilde{N}$.
By (\ref{eq:keyes}) we have $p^{d}\tilde{N} \subseteq \iota(N)$, so that
this implies
that $p^{a}x - p^{a+b}n \in \ker \iota$ for some $n \in N$.
Again by (\ref{eq:keyes}) this kernel is killed by $p^{d}$; we conclude that
$$p^{a+d}x \in p^{a+b+d}N.$$
Thus
$$x \in p^{b}N + N[p^{a+d}] \subseteq p^{b}N \subseteq \p_{b}N.$$
Since $x \notin \p_{b}N$ by hypothesis, this yields the desired contradiction.
\end{proof}

\section{Reductions of Mordell-Weil groups} \label{s:aa}

\subsection{Galois elements} \label{s:a3}

Let $A$ be an abelian variety over a number field $F$.  By
\cite[Section 19, Corollary 2]{Mumford} the ring
$\O := \Ef $ is a reduced, finite, flat $\Z$-algebra.  We further assume
that it is commutative; we fix a rational prime $p$, and we continue with the
notations of Section~\ref{s:a} for this ring $\O$ and prime $p$.
By (\ref{eq:l1}) we may
fix $a_{i,n} \in \p_{i,n}$ and $b_{i,n} \in \prod_{j \neq i} \p_{j,n}$ 
such that $a_{i,n} + b_{i,n} = c^{g-1}p^{d(g-1)}$.  The map
\begin{align*}
\varphi_{n} : A[p^{n-d}] &\to A[\p_{1,n}] \oplus \cdots \oplus A[\p_{g,n}] \\
t &\mapsto (b_{1,n}t,\ldots,b_{g,n}t)
\end{align*}
is then well-defined by (\ref{eq:l2}).

\begin{lemma} \label{l:coker}
The cokernel of $\varphi_{n}$ is bounded independent of $n$.
\end{lemma}
\begin{proof}
Since $p^{n} \in \p_{i,n}$ we can define a map
\begin{align*}
\psi_{n} : A[\p_{1,n}] \oplus \cdots \oplus A[\p_{g,n}] &\to A[p^{n-d}] \\
(t_{1},\ldots,t_{g}) &\mapsto p^{d}(t_{1}+\cdots+t_{g}).
\end{align*}
As $c^{g-1}p^{d(g-1)} - b_{i,n} \in \p_{i,n}$, the map
$\varphi_{n} \circ \psi_{n}$ is just multiplication
by $c^{g-1}p^{dg}$.  The lemma follows from this.
\end{proof}

For an $\O$-submodule $\Gamma$ of $A(F)$, we now write
$$\lambda_{\p_{i,n+d}}^{\Gamma} : \gpns \to
\Hom_{\O}\bigl(\Gamma,A[\p_{i,n+d}]\bigr)$$
for the composition of $\lpns$ with $\varphi_{n+d}$ and projection to
$A[\p_{i,n+d}]$.
In the next lemma we use the natural map $\gpns \to \g_{p^{m}}(\Gamma)$ 
(corresponding to multiplication
by $p^{n-m}$ from $\Hom_{\O}(\Gamma,A[\p_{i,n+d}])$ to
$\Hom_{\O}(\Gamma,A[\p_{i,m+d}])$)
to regard $\lambda_{\p_{i,m+d}}^{\Gamma}$ as a map from
$\gpns$ for $n \geq m$.

\begin{lemma} \label{l:galois}
Let $y_{1},\ldots,y_{r}$ be an $\O$-pre-basis of $A(F)$.
Then there is an integer $b$ such that for all sufficiently large $n$ there
is a $\sigma_{n} \in \gpn(A(F))$ with
$$\lpn^{\O \cdot y_{j}}(\sigma_{n}) = 0 \text{~for~} j=2,\ldots,r;$$
$$\lambda_{\p_{i,b}}^{\O \cdot y_{1}}(\sigma_{n}) \neq 0
\text{~for all~} i \in I_{p}(y_{1}).$$
\end{lemma}
\begin{proof}
The cokernel of the natural map
$$\pi : \Hom_{\O}\bigl(A(F),A[p^{n}]\bigr) \to \op{j=1}{r}
\Hom_{\O}\bigl(\O \cd y_{j},A[p^{n}]\bigr)$$
is bounded independent of $n$ by the definition of a pre-basis.
Combined with Proposition~\ref{p:padic}, it follows that the cokernel of
$$\pi \circ \lpn^{A(F)} : \gpn\bigl(A(F)\bigr) \to
\op{j=1}{r} \Hom_{\O}\bigl(\O \cd y_{j},A[p^{n}]\bigr)$$
is bounded independent of $n$.  Finally, by Lemma~\ref{l:coker} we conclude
that the cokernel of the map
\begin{multline} \label{eq:big}
\gpn\bigl(A(F)\bigr) \to \\
\Bigl(\underset{i \in I_{p}(y_{1})}
{\oplus}\Hom_{\O}\bigl(\O \cd y_{1},
A[\p_{i,n+d}]\bigr) \Bigr) 
\oplus \Bigl( \op{j=2}{r} \Hom_{\O}\bigl(\O \cd y_{j},A[p^{n}]
\bigr) \Bigr)
\end{multline}
is bounded independent of $n$.  

By the definition of the set $I_{p}(y_{i})$, 
for each $i \in I_{p}(y_{1})$ there is some $m > 0$
such that $p^{n+d-m}\Hom_{\O}(\O \cd y_{1},A[\p_{i,n+d}]) \neq 0$ for
sufficiently large $n$.  (That is, these groups grow with $n$.)
Since the cokernel of (\ref{eq:big}) is bounded,
it follows that there is an integer
$b$ such that for sufficiently large $n$ there is $\sigma_{n} \in
\gpn(A(F))$ with
$$\sigma_{n}|_{\Hom_{\O}(\O \cd y_{j},A[p^{n}])} = 0 \text{~for~} 
j=2,\ldots,r;$$
$$p^{n+d-b}\sigma_{n}|_{\Hom_{\O}(\O \cd y_{1},A[p_{i,n+d}])} \neq 0
\text{~for all~} i \in I_{p}(y_{1}).$$
By the remarks preceding the lemma, this $\sigma_{n}$ 
is the required element of $\gpn(A(F))$.
\end{proof}

\begin{lemma} \label{l:place}
Let $y_{1},\ldots,y_{r}$ be an $\O$-pre-basis of $A(F)$.  Then there is an
integer $b$ such that for all sufficiently large $n$ there are infinitely many
places $w$ of $F_{p^{n}}$ with
$$\rw  y_{j} \in p^{n}A(k_{w}) \text{~for~} j = 2,\ldots,r;$$
$$\rw y_{1} \notin \p_{i,b}A(k_{w}) \text{~for~} i \in I_{p}(y_{1}).$$
\end{lemma}
\begin{proof}
Let $n$ be sufficiently large and fix
$\sigma_{n}$ as in Lemma~\ref{l:galois}.  If $w$ is a place of
$F_{p^{n}}$ with $\Frob_{w} = \sigma_{n}$ in $\g_{p^{n}}(A(F))$, then $w$
satisfies the conditions of the lemma
by Lemmas~\ref{l:frob} and~\ref{l:crt2}.
Since the Cebatorev density theorem guarantees the existence of infinitely
many such $w$, the lemma follows.
\end{proof}

\subsection{Reduction of subgroups} \label{s:a4}

We are now in a position to prove our main result.

\begin{proposition} \label{p:main}
Let $A$ be an abelian variety over a number field $F$; 
assume that $\O = \Ef A$ is commutative.  Fix a rational prime $p$ and
let $\Sigma$ be a subgroup of $A(F)$ containing $A(F)_{\tors}$.
Suppose that $x \in A(F)$ is such that
\begin{equation} \label{eq:red}
\rv  x \in \rv  \Sigma
\end{equation}
for almost all places $v$ of $F$.  Then $x$ lies in $\Sigma \otimes \Zpp$.
\end{proposition}
\begin{proof}
Suppose that $x \notin \Sigma \otimes \Zpp$.  By
Proposition~\ref{p:pb} we can then
choose an  $\O$-pre-basis $y_{1},\ldots,y_{r}$ of
$A(F)$ such that there is an integer $a$ with
\begin{equation} \label{eq:lat}
\psi_{1}(x) \notin \psi_{1}(\Sigma) + p^{a}\O.
\end{equation}

Let $b$ be the integer determined by $y_{1},\ldots,y_{r}$ in
Lemma~\ref{l:place} and fix $n > a + b + 2d$.  Let $w$ be a place
of $F_{p^{n}}$ as in Lemma~\ref{l:place};
by (\ref{eq:red}) we may further assume that there is a $y \in \Sigma$ with
$\rw  x = \rw  y$.  
Multiplying by $\eta$, by (\ref{eq:eta}) we have
$$\psi_{1}(x)\rw y_{1} + \cdots + \psi_{r}(x)\rw y_{r} =
\psi_{1}(y)\rw y_{1} + \cdots + \psi_{r}(y)\rw y_{r}.$$
Thus
\begin{equation} \label{eq:red3}
\bigl(\psi_{1}(x) - \psi_{1}(y) \bigr) \rw y_{1} \in p^{n}A(k_{w})
\end{equation}
by the definition of $w$.

Set $\alpha = \psi_{1}(x) - \psi_{1}(y)$; by (\ref{eq:lat}) and
(\ref{eq:l3}), $\alpha \notin \p_{i,a+d}$ for some $i$.
Fix such an $i$.  Since $\alpha \in \im \psi_{1}$, by (\ref{eq:pbind}) we have
$i \in I_{p}(y_{1})$; thus we also have
$\rw y_{1} \notin \p_{i,b}A(k_{w})$ by
the definition of $w$.
Since $A(k_{w})[p^{a+2d}] \subseteq p^{b}A(k_{w})$ (as
$A[p^{n}] \subseteq A(k_{w})$ and $a+b+2d < n$), we may therefore
apply Lemma~\ref{l:evil} to conclude that
$\alpha x \notin p^{a+b+2d}A(k_{w})$.
Since $a + b + 2d < n$, this contradicts (\ref{eq:red3}), and thus proves
the proposition.
\end{proof}

\begin{corollary}
Let $A$ be an abelian variety over a number field $F$ and assume that
$\Ef A$ is commutative.  Let $\Sigma$ be a subgroup of $A(F)$
containing $A(F)_{\tors}$ and suppose that $x \in A(F)$ is such that
$\rv x \in \rv \Sigma$
for almost all places $v$ of $\Sigma$.  Then $x \in \Sigma$.
\end{corollary}
\begin{proof}
This is immediate from Proposition~\ref{p:main} applied for all primes
$p$.
\end{proof}

\end{document}